\documentclass[noinfoline]{imsart}

\RequirePackage{amsthm,amsmath}
\RequirePackage[numbers]{natbib}
\RequirePackage[colorlinks,citecolor=blue,urlcolor=blue]{hyperref}
\usepackage{imsart}
\usepackage{amscd}
\usepackage{mathrsfs}
\usepackage{latexsym}
\usepackage{stmaryrd}
\usepackage{xypic}
\xyoption{curve}
\usepackage{wasysym}
\usepackage{enumitem}
\usepackage{verbatim}
\usepackage{ifthen}
\usepackage{amsfonts,amssymb,booktabs,color,epsfig,graphicx,hyperref,nicefrac,url,array,rotating}
\usepackage[margin=1.5in]{geometry}

\newtheorem{thm}{Theorem}
\newtheorem{cor}{Corollary}
\newtheorem{lem}{Lemma}

\theoremstyle{definition}

\newtheorem{eg}{Example}

\newcommand{\SL}{\mathbb{SL}} 
\newcommand{\SO}{\mathbb{SO}} 
\newcommand{\X}{\mathbf{X}} 
\newcommand{\G}{\mathbf{G}} 
\newcommand{\K}{\mathbf{H}} 
\newcommand{\A}{\mathbf{A}} 
\newcommand{\B}{\mathbf{B}} 
\newcommand{\N}{\mathbf{N}} 
\newcommand{\M}{\mathbf{M}} 
 
\newcommand{\R}{\mathbb{R}} 
\newcommand{\C}{\mathbb{C}} 
\newcommand{\HH}{\mathbb{H}_2} 
\newcommand{\half}{\nicefrac{1}{2}} 
\newcommand{\specmeasure}{\tau}
\newcommand{\empchar}{\hat\phi} 
\newcommand{\technical}[1]{#1} 

\begin{document}
\begin{frontmatter}
\title{Kernel Density Estimation on Symmetric Spaces of Non-Compact Type}
\runtitle{KDE for Symmetric Spaces}

\begin{aug}
\author{\fnms{Dena Marie} \snm{Asta} }
	\address{Department of Statistics \\ The Ohio State University\\ 1958 Neil Ave\\ Columbus, OH 43210 USA\\ \href{Url}{dasta@stat.osu.edu}}
\runauthor{D. M. Asta}
\end{aug}

\begin{abstract}
  We construct a kernel density estimator on symmetric spaces of non-compact type and establish an upper bound for its convergence rate, analogous to the minimax rate for classical kernel density estimators on Euclidean space.  
  Symmetric spaces of non-compact type include hyperboloids of constant curvature $-1$ and spaces of symmetric positive definite matrices.
  This paper obtains a simplified formula in the special case when the symmetric space is the space of normal distributions, a 2-dimensional hyperboloid.
\end{abstract}


\tableofcontents

\end{frontmatter}

\section{Introduction}

Data, while often expressed as collections of real numbers, are often more naturally regarded as points in non-Euclidean spaces.
To take an example, radar systems can yield the data of bearings for planes and other flying objects; those bearings are naturally regarded as points on a sphere \citep{Rahman-et-al}.
To take another example, diffusion tensor imaging (DTI) can yield information about how liquid flows through a region of the body being imaged; that three-dimensional movement can be expressed in the form of symmetric positive definite $(3\times 3)$-matrices \citep{Rahman-et-al}.
To take yet another example, the nodes of certain hierarchical real-world networks can be regarded as having latent coordinates in a hyperboloid \citep{Krioukov-et-al-hyperbolic-geometry, Asta-et-al}.
In all such examples, the spaces can be regarded as subsets of Euclidean space even though Euclidean distances do not reflect true distances between points.
An ordinary kernel density estimator (KDE) applied to sample data generally will not be optimal in terms of the $L_2$-risk with respect to the volume measure on the non-Euclidean manifold.

The idea of kernel density estimation is to smooth out, or convolve, an empirical estimator (an average of dirac distributions centered at the data) with a smooth rapidly decaying \technical{kernel} so as to obtain a smooth estimate of the true density.
The literature offers some variants of kernel density estimation on non-Euclidean spaces.
A simple variant, for compact manifolds \citep{Pelletier} or more general compact subsets of manifolds \citep*{chevallier2016kernel,chevallier2017kernel}, applies a Euclidean kernel having supports small enough to fit inside the charts.
A more general version, defined on complete manifolds, generalizes the kernel to be defined on the tangent bundle - in effect, requiring a kernel for each point \citep{Kim-et-al}.
Minimax rates of convergence have been proven for all of these different variants in terms of a H\"{o}lder class exponent \citep{Pelletier,Kim-et-al}.
It is desirable to refine these convergence rates based on Sobolev constraints on the true densities.

On \technical{symmetric spaces} like Euclidean space, a kernel need only be defined at one point and transported everywhere else.
All symmetric spaces $X$ can be decomposed into symmetric spaces of \technical{Euclidean}, \technical{compact}, and \technical{noncompact} type in such a way that a KDE on $X$ can be constructed from KDEs on the three types.   
Symmetric spaces of the first two type admit KDEs with minimax convergence rates.
The goal of this paper is to begin to complete the picture for symmetric spaces of noncompact type, by proving the upper bound part of a conjectured minimax rate.  



Kernel density estimation for random variables can be interpreted as an estimate of a Fourier transform of the density.
Thus for manifolds on which Fourier analysis generalizes, there should exist some generalization of the KDE.
One of the earliest such Fourier-based generalizations of the KDE is defined for compact manifolds and shown to be minimax \citep{Hendriks}. 
Certain density estimators on the non-compact \technical{Poincar\'{e} halfplane} \citep{Huckemann-et-al} and the space of symmetric positive definite matrices, based on \technical{Helgason-Fourier Analysis}, have been shown to be minimax for estimation from corrupted samples.
In the special case where the noise is non-existent, these estimators can be regarded as special cases of a KDE where the kernel is a natural generalization of the sinc kernel.
However, these estimators have not been shown to be minimax for estimation from uncorrupted samples.
Moreover, the full of generality of Helgason-Fourier Analysis is not exploited in defining and analyzing this estimator.

The Helgason-Fourier transform, unlike the ordinary Fourier transform, sends functions on a symmetric space of noncompact type to functions on a different frequency space.
The exact form of this frequency space, much less a usable formula for the transform, depends on a geometric understanding the original space.
Countless symmetric spaces of interest in applications have well-understood geometries.
When the original space is the Poincar\'{e} halfplane, for example, the frequency space is a cylinder.
The Helgason-Fourier transform, like the ordinary Fourier transform, is an isometry on $L_2$-function spaces and sends convolutions to products in a certain sense.

This paper uses the Helgason-Fourier transform to construct and analyze a version of a KDE on symmetric spaces of noncompact type.
We define a kernel density estimator for symmetric spaces of non-compact type $\X$, for which Helgason-Fourier transforms are defined.
Unlike the non-Fourier-based variant \citep{Pelletier} for compact manifolds, this variant is differentiable everywhere and estimates densities with non-compact support.
The analogue of a kernel is often just a density on a space $\G$ of isometries invariant with respect to the subgroup $\K$ of $\G$ for which $\X=\G/\K$.   
An example is a Gaussian, a solution to the $\K$-invariant heat equation,
We bound risk in terms of bandwidth $h$, the number $n$ of sample points, a Sobolev parameter $\alpha$, and the sum of the restricted roots of $\X$ [Theorem \ref{thm:upper}].
Optimizing $h$ in terms of $n$, we obtain an upper bound of
$$n^{-2\alpha/(2\alpha+\dim{\X})},$$
for the convergence rate, where $\alpha$ is a Sobolev parameter, under natural assumptions on the density space [Theorem \ref{thm:upper}].
We then obtain a simplified formula, that can be implemented on a computer, for the special case where $\X$ is the $2$-dimensional hyperboloid of constant curvature $-1$.
The proof for the upper bound of the convergence rates adapts techniques used in \citep{Huckemann-et-al}. 
We conjecture that the same upper bound yields a lower bound and hence a minimax rate, and a proof is reserved for future work.

\section{Preliminaries}
We recall some preliminary constructions and results in this section.
First, we recall the construction of a KDE, including a Fourier-based interpretation.
We then recall some of the theory of symmetric spaces.
Finally, we recall the theory of Helgason-Fourier Analysis.

\subsection{Kernel Density Estimation}
The \technical{kernel density estimator} (KDE) $f^{(h)}_{(X_1,\ldots,X_n)}:\R\rightarrow\R$, defined by
\begin{equation}
  \label{eqn:euclidean.kde}
  f^{(h)}_{(X_1,X_2,\ldots,X_n)}(x)=\frac{1}{nh}\sum_{i=1}^nK\left(\frac{x-X_i}{h}\right),
\end{equation}
estimates a density $f$ on $\R$ based on some observed points $X_1,\ldots,X_n$ identically and independently sampled from $f$, a tunable \technical{bandwidth} $h>0$, and a fixed choice of \technical{kernel} function $K$, which in most cases amounts to a symmetric unimodal Lebesgue density on $\R$ having mode $0$.

A generalization of the KDE for certain Riemannian manifolds $X$ is
\begin{equation}
  \label{eqn:simple.manifold.kde}
  f^{(h)}_{(X_1,X_2,\ldots,X_n)}(x)=\frac{1}{nh^{\dim X}}\sum_{i=1}^n\theta_{X_i}(x)^{-1}K\left(\frac{\mathrm{dist}(x-X_i)}{h}\right),
\end{equation}
where $\theta_p$ denotes the density of the volume measure on the Riemannian manifold \citep{Pelletier}.
For Euclidean space, $\theta_p$ is just the constant $1$ function and $K$ can be taken to be a general kernel function.
For compact manifolds, $K$ is taken to have support $[-1,+1]$ and bandwidth $h$ is bounded by the injectivity radius of the manifold.
In these cases, (\ref{eqn:simple.manifold.kde}) integrates to $1$, defines a density when $K$ is non-negative, and converges to $f$ at a minimax rate in a suitable sense \citep{Pelletier}.

One way to think about kernel density estimation is that it is an estimation of the \technical{characteristic function} $\phi_X$ of a random variable $X$, defined by
\begin{equation}
  \label{intro:eqn:characteristic}
  \phi_X(t)=\mathbb{E}[e^{itX}].
\end{equation}
An \technical{empirical characteristic function} $\empchar_{x_1,\ldots,x_n}$, defined by
\begin{equation}
  \label{intro:eqn:empirical.characteristic}
  \empchar_{x_1,\ldots,x_n}(t)=\frac{1}{n}\sum_{i=1}^ne^{itx_i},
\end{equation}
estimates $\phi_X$ based on some observed points $x_1,\ldots,x_n$ sampled from $X$.
While (\ref{intro:eqn:characteristic}) often admits a convergent inverse Fourier transform, (\ref{intro:eqn:empirical.characteristic}) does not.
Therefore while the density for $X$ can often be recovered by taking the inverse Fourier transform of (\ref{intro:eqn:characteristic}), an estimated characteristic function for $X$ does not analogously give an estimated density for $X$.
Instead, (\ref{intro:eqn:empirical.characteristic}) needs to be dampened by a rapidly decaying function $\mathcal{F}[K_h]$, after which it is in $L_2(\R)$ and therefore has a well-defined inverse Fourier transform. 
If $\mathcal{F}[K]$ is fixed and $\mathcal{F}[K_h](s)$ is set to be $\mathcal{F}[K](hs)$, one then recovers the original definition (\ref{eqn:euclidean.kde}) of a KDE.
What is more, Fourier Analysis then makes it then possible to prove that the standard KDE achieves a minimax rate for densities in a Sobolev ball.
Thus different types of Fourier transforms give different variants of the standard KDE.  
In this way, an alternative to (\ref{eqn:simple.manifold.kde}) for certain compact spaces achieves a minimax rate for densities in a Sobolev ball.

\subsection{Symmetric spaces}
This paper assumes the basic definition of a \technical{smooth manifold}.  
A Riemannian manifold is a smooth manifold $M$ equipped with a \technical{Riemannian metric}, a choice of inner product $g_x:T_xM\otimes T_xM\rightarrow T_xM$ on the real tangent vector spaces $T_xM$ smoothly varying in $x\in M$.  
The formula  $\int_t\langle\gamma'(t),\gamma'(t)\rangle\;dt$ for arc lengths of curves $\gamma:[0,1]\rightarrow\R^n$ straightforwardly generalizes to a definition of arg lengths of curves, and hence of geodesics and distances between points, in Riemannian manifolds.  
An \technical{isometry} between Riemannian manifolds is a smooth map of manifolds whose differential defines an isometry of tangent spaces.
An \technical{involution} at a point $x$ in a Riemannian manifold $M$ is an isometry $\varphi:M\cong M$ fixing $x$ such that the induced linear isometry $(\partial\varphi)_x:T_xM\rightarrow T_xM$ is defined by scalar multiplication by $-1$.  
A \technical{Riemannian symmetric space} is a Riemmaniain manifold admitting an involution at each of its points.
For background on the basic theory of Riemannian symmetric spaces, the reader is referred to \cite{loos1969symmetric}.  

If a Riemannian symmetric space $X$ decomposes as a product $X=Y\times Z$ of such spaces, then KDEs on the factors $Y$ and $Z$ can be combined to give a KDE on the entire space.
If a Riemmanian symmetric space $X$ is a quotient $E/G$ of a Riemannian symmetric space $E$ by the action of a discrete group $G$, then a KDE on $X$ can be defined as the restriction of a KDE on $E$ by identifying all but a Borel measure $0$ subspace of $X$ with a subspace of $E$.  
Therefore it suffices to restrict our attention to \technical{irreducible simply connected Riemannian symmetric spaces}, Riemannian symmetric spaces that do not admit decompositions into a product of non-trivial Riemannian symmetric spaces and do not arise as quotients $E/G$ of Riemannian symmetric spaces $E$ by non-trivial discrete groups $G$.  
The irreducible simply connected Riemannian symmetric spaces fall into three types: Euclidean, compact, and non-compact (and non-Euclidean).  
A minimax KDE on Euclidean space is classical.
A minimax KDE on compact Riemannian manifolds already exists.
This paper focuses on the (irreducible, simply connected, and) noncompact case.

Such spaces admit the following algebraic characterization.
Recall that a \technical{Lie group} is a smooth manifold $G$ that is at once a group in such a way that the  multiplication $G\times G\rightarrow G$ and inversion $G\rightarrow G$ are smooth maps.  
Let $\G$ be a noncompact \technical{semisimple Lie group}, a noncompact connected Lie group with finite center containing no non-trivial normal connected Abelian subgroups.
Then $\G$ admits an \technical{Iwasawa decomposition} $\G=\K\A\N$, a decomposition of $\G$ as the space $\K\A\N$ of all triple products of elements from a maximal compact Lie subgroup $\K\leqslant\G$, an Abelian Lie subgroup $\A\leqslant\G$, and  \technical{nilpotent} Lie subgroup $\N\leqslant\G$.
(Since the symbol $K$ will assume its traditional role in statistics as denoting a kernel, the symbol $\K$ is used here to denote a maximal compact subgroup of $\G$.)
Nilpotency means that $\N_k$ is the trivial group for large enough $k\gg 0$, where $\N_0=\N$ and $\N_{i+1}$ is the subgroup of $\N$ generated by all elements of the form $xyx^{-1}y^{-1}$ for $x\in\N_i$ and $y\in\N$.
In practice, $\G$ is a suitable group of matrices under matrix multiplication, $\A$ is a group of diagonal matrices and $\N$ is a group of upper triangular matrices.  

\begin{eg}
  \label{eg:sl2-iwasawa}
  The group $\SL_2$ has Iwasawa decomposition
  $$\SL_2=\SO_2\mathbb{R}_+\mathbb{R},$$
  where the positive reals $\mathbb{R}_+$ (under the operation of multiplication) is identified with the diagonal $(2\times 2)$-matrices in $\SL_2$ with positive entries and the reals $\mathbb{R}$ (under the operation of addition) is identified with the upper triangular $(2\times 2)$-matrices with $1$'s along the diagonal.
\end{eg}

Let $\X$ be the smooth manifold defined as the quotient space
$$\X=\G/\K\cong\A\N.$$
The space $\X$ is a smooth manifold with smooth structure characterized by the property that precomposition with the natural function $\G\rightarrow \G/\K$ bijectively identifies smooth functions $\G/\K\rightarrow\mathbb{R}$ with smooth functions $\phi:\G\rightarrow\mathbb{R}$ which are \technical{right-$\K$-invariant} ($\phi(gh)=\phi(g)$ for all $h\in\K$).
A choice of a bi-$\K$-invariant, left $\G$-invariant inner product on $\mathfrak{g}$ passes to a well-defined $\G$-invariant inner product on a tangent space of $\X$, which in turn uniquely extends to a Riemannian metric turning $\X$ into a Riemannian symmetric space.  
Conversely, every Riemannian symmetric space of noncompact type arises in this manner.

\begin{eg}
  \label{eg:H2}
  Continuing Example \ref{eg:sl2-iwasawa}, we can give an algebraic construction of the \technical{Poincar\'{e} halfplane} $\HH$, the $2$-manifold defined as the subspace
  $$\HH=\{z\in\mathbb{C}\;|\;\mathrm{Im}(z)>0\},$$
  of $\C$ equipped with the Riemannian metric given by the arc length
  $$ds^2=(\mathrm{Im}\;z)^{-2}(d(\mathrm{Re}\;z)^2+d(\mathrm{Im}\;z)^2).$$
  The matrices in $\SL_2$ act on $\HH$ as \technical{M\"{o}bius transformations}:
  $$\left( \begin{array}{cc} a & b \\ c & d \end{array} \right)(z)=\frac{az+b}{cz+d}.$$
  The matrices in $\SL_2$ fixing $i\in\HH$ form the matrix subgroup $\SO_2$.
  The action of $\SL_2$ on $\HH$ implicitly gives a well-defined bijection $\SL_2/\SO_2\cong\HH$ sending an equivalence class of a matrix $m\in\SL_2$ to $m(i)\in\HH$.
  This bijection $\HH\cong\SL_2/\SO_2$ defines an isometry for a suitable choice of bi-$\SO_2$-invariant inner product on the Lie algebra $\mathfrak{sl}_2$ associated to $\SL_2$ \cite[\S3.1]{Terras-harmonic-1}.
  Thus $\HH$ is a Riemmanian symmetric space.  
  This space can be interpreted as the information manifold of all univariate normal distributions, where the real coordinates describe the means, and the imaginary coordinates describe standard deviations, and the Riemannian metric is the Fisher metric.  
  Alternatively, this space is a natural latent space for families of random graphs used to model real-world networks \cite{Krioukov-et-al-hyperbolic-geometry}. 
  Alternatively, this space models electrical impedances on which certain circuit elements act as M\"{o}bius transformations \citep{Huckemann-et-al}.
\end{eg}

From now on, fix a semisimple Lie group $\G$ having finite center with Iwasawa decomposition $\G=\K\A\N$ and let $\X=\G/\K$.
Let $\mathfrak{g},\mathfrak{a},\mathfrak{n}$ denote the Lie algebras of the respective Lie groups $\G,\A,\N$, their tangent spaces over identities.  
Recall that the \technical{exponential map} $\exp:\mathfrak{a}\rightarrow\A$ is defined as sending a tangent vector $v$ to $\gamma_v(1)\in\A$ for $\gamma_v$ the unique continuous group homomorphism $\gamma:\R\rightarrow\A$ with $\gamma'(0)=v$.
For each $x\in\X$, there exist unique elements $a(x)\in\mathfrak{a}$ and $n(x)\in\N$ such that $x=\exp(a(x))n(x)$.

\begin{eg}
  \label{eg:H2.coordinates}
  Continuing Examples \ref{eg:sl2-iwasawa} and \ref{eg:H2}, for each $z\in\HH$, 	
	$$n(z)=\left(\begin{tabular}{cc}$1$&$\mathrm{Re}\;z$\\$0$&$1$\end{tabular}\right),\quad\exp(a(z))=\left(\begin{tabular}{cc}$\sqrt{\mathrm{Im}\;z}$&$0$\\$0$&$(\sqrt{\mathrm{Im}\;z})^{-1}$\end{tabular}\right).$$
\end{eg}

Let $[-,-]:\mathfrak{g}\times\mathfrak{g}\rightarrow\mathfrak{g}$ denote the \technical{Lie bracket} operation, defined in the case $\G$ is a group of $(n\times n)$-invertible matrices and thus $\mathfrak{g}$ is a vector space of $(n\times n)$-matrices, by $[v,w]=vw-wv$.  
Recall that the \technical{Killing form} on $\mathfrak{g}$ is the symmetric bilinear map $\kappa:\mathfrak{g}\times\mathfrak{g}\rightarrow\R$ sending a pair $(v,w)$ of tangent vectors to the trace of the operator $[x,[y,-]]$ on $\mathfrak{g}$.
The Killing form on $\mathfrak{g}$ restricts to an inner product on $\mathfrak{a}$.  
In this manner, $\mathfrak{a}^*$ will sometimes be naturally identified with $\mathfrak{a}$ along the adjoint of this inner product, and $\mathfrak{a}$ will often by identified with the Lebesgue measure space $\R^{\dim\mathfrak{a}}$ along an isometry between the two.  
Let $\M$ be the set of all elements in $\K$ commuting with all elements in $\A$.
Then let $\B=\K/\M$.
For each $\lambda\in\mathfrak{a}^*$, let $\mathfrak{g}_\lambda$ denote the space of all $g\in\mathfrak{g}$ with $\lambda(x)g=xg-gx$ for all $x\in\mathfrak{a}$.
The \technical{restricted root system} $\Lambda$ of $\G$ is the set of all $\lambda\in\mathfrak{a}^*$ with $\dim\mathfrak{g}_\lambda>0$.
A \technical{positive restricted root system} is a choice of subset $\Lambda_+\subset\Lambda$ such that $\Lambda=\Lambda^+\cup-\Lambda^+$ and $\Lambda^+\cap-\Lambda^+=\varnothing$.  
Set
\begin{equation*}
	w_{\X}=\#\{k\M\;|\;k\A=\A k,\;k\in\K\},\quad\rho_{\X}=\frac{1}{2}\sum_{\lambda\in\Lambda^+}\dim\mathfrak{g}_\lambda.
\end{equation*}
The root systems $\Lambda$ play a pivotal role in the classification of symmetric spaces.
Note that the choice of $\Lambda_+$, and hence also of $\rho_{\X}$, is only unique up to sign.

\begin{eg}
	For $\mathbb{H}_2=\SL_2/\SO_2$, we have the following \cite[\S3.1]{Terras-harmonic-1}:
	$$\rho_{\HH}=\half,\;w_{\mathbb{H}_2}=1.$$ 
\end{eg}

\subsection{Helgason-Fourier Analysis}
Helgason-Fourier Analysis is an analogue of Fourier Analaysis for symmetric spaces of noncompact type.  
The reader is referred to \citep[Section 2]{Pesenson} for a concise summary of the theory and \citep{Terras-harmonic-1} for details in the special case $\X=\HH$.  
The \technical{Helgason-Fourier transform} of $f\in L_1(\X,dx)$ at $\lambda\in\mathfrak{a}^*$ and $b=k\M\in\B=\K/\M$, when it exists, is
$$(\mathcal{H}{f})(\lambda,b)=\int_{x=g\K\in\X=\G/\K}f(x)e^{{\overline{-s(a(k^{-1}g))}}}\;dx,$$
where $s=\rho_{\X}-i\lambda$ lies in the complexification of $\mathfrak{a}^*$, and $dx$ is the intrinsic volume measure on $X$.
The Helgason-Fourier transform defines a linear map
\begin{equation*}
  \mathcal{H}:L_2({\X},dx)\rightarrow L_2\left(\mathfrak{a}^*\times\B,|c|^{-2}d\lambda\;db\right).
\end{equation*}
such that we have the following Plancharel identity
\begin{equation}
  \label{plancherel}
  \int_{\X}|f(x)|^{2}dx=w_{\X}^{-1}\int_{\lambda\in\mathfrak{a}^{*}}\int_{b\in\B}\left|\mathcal{H}f\left(\lambda,b\right)\right|^{2}|c(\lambda)^{-2}|\;d\lambda\;db.
\end{equation}

Here $d\lambda$ is the Lebesgue measure under the natural identification $\mathfrak{a}^*\cong\mathbb{R}^r$, $db$ is a suitably normalized $K$-invariant measure, and $c$ is the \technical{Harish-Chandra c-function} $\mathfrak{a}^*\rightarrow\R$ \citep{helgason1994harish}; a proper treatment of the latter, while fundamental throughout geometric analysis, is beyond the scope of the paper.
However, the following exact formula is used in the special case $\X=\HH$.

\begin{eg}
  Continuing Examples \ref{eg:H2} and \ref{eg:sl2-iwasawa},
  $$c(\lambda)=\frac{1}{8\pi^2}\lambda\tanh (\pi \lambda).$$
\end{eg}

The basic property of the Harish-Chandra c-function needed in the proofs is the following growth bound, noted in \citep{Narayanan-et-al}.

\begin{lem}
  \label{lem:c-bound}
  For each $\lambda\in\mathfrak{a}^*$, $|c(\lambda)|^{-2}\leqslant(1+|\lambda|)^{\dim\mathfrak{n}}$.
\end{lem}

The Helgason-Fourier transform $\mathcal{H}$ sends convolutions to products:
\begin{equation*}
  \mathcal{H}[f_{GX}]=\mathcal{H}[f_G]\mathcal{H}[f_X]
\end{equation*}
for each bi-$\K$-invariant density $f_{G}\in L_2(\G)$ of a random quantity $G$ on $\G$, identified with the induced left-$\K$-invariant density on $\X$, and density $f_X$ of a random quantity $X$ on $\X$, where $f_{GX}$ is the density associated to $GX$.  

The Sobolev ball $\mathcal{F}_\alpha(Q)$ in $L_2(\X)$ can be defined as
$$\mathcal{F}_\alpha(Q)=\{f\in L_2(\X)\;|\;\|\Delta^{\alpha/2}f\|^2\leqslant Q\},$$
just as in the Euclidean case, except that here $\Delta^{\alpha/2}f$ is defined by
$$\mathcal{H}[\Delta^{\alpha/2}f](\lambda,b)=(-\rho_{\X}^2-\lambda^2)^{\alpha/2}\mathcal{H}[f_X](\lambda,b)$$

The \technical{inverse Helgason-Fourier transform}, written $\mathcal{H}^{-1}$, is given by
$$\label{inverse}\mathcal{H}^{-1}\psi(x)=w_{\X}^{-1}\int_{\mathfrak{a}^{*}}\int_{\B}\psi(\lambda,k\M)e^{{(\rho_{\X}-i\lambda)(a(k^{-1}x))}}\;|c|^{-2}d\lambda\;db,$$

We can recover $f\in\mathbb{C}_{c}^{\infty}(\X)$ from its transform by the inversion formula
$$f=\mathcal{H}^{-1}\mathcal{H}f.$$

\subsection{$\G$-kernel density estimation}
In an analogy with the classical definition, define the \technical{characteristic function}
$$\phi_X:\mathfrak{a}^*\times\B\rightarrow\C$$ 
of a random quantity $g\K$ on the symmetric space $\X=\G/\K$ by
\begin{equation*}
  \phi_X(s,k\M)=\mathbb{E}[e^{\bar{s}(-a(k^{-1}g))}].
\end{equation*}
An \technical{empirical characteristic function} $\empchar_{x_1,\ldots,x_n}:\mathfrak{a}^*\rightarrow\C$, defined by
\begin{equation}
  \label{eqn:empirical.characteristic}
  \empchar_{g_1\K,\ldots,g_n\K}(s,k\M)=\frac{1}{n}\sum_{i=1}^ne^{\bar{s}(-a(k^{-1}g_i))}.
\end{equation}
estimates $\phi_X$ based on some observed points $g_1\K,\ldots,g_n\K$ sampled from a random quantity on $\X$.
This function needs to be dampened by a \technical{kernel} before taking $\mathcal{H}^{-1}$.
A \technical{$\G$-kernel $K$ on $\X$} is a unimodal left-$\K$-invariant $L_2$-function $K$ on $\X$ which integrates to $1$.

A non-density example is an analogue $K_{\mathrm{sinc}}$ of the sinc kernel defined by
\begin{equation*}
  \hat{K}_{\mathrm{sinc}}(s,b)\propto
  \begin{cases}
		1 & |s|<1,\\
    0 & |s|>1.
  \end{cases}
\end{equation*}

\begin{figure}[h]
  \begin{center}
		{\includegraphics[width=60mm,height=50mm]{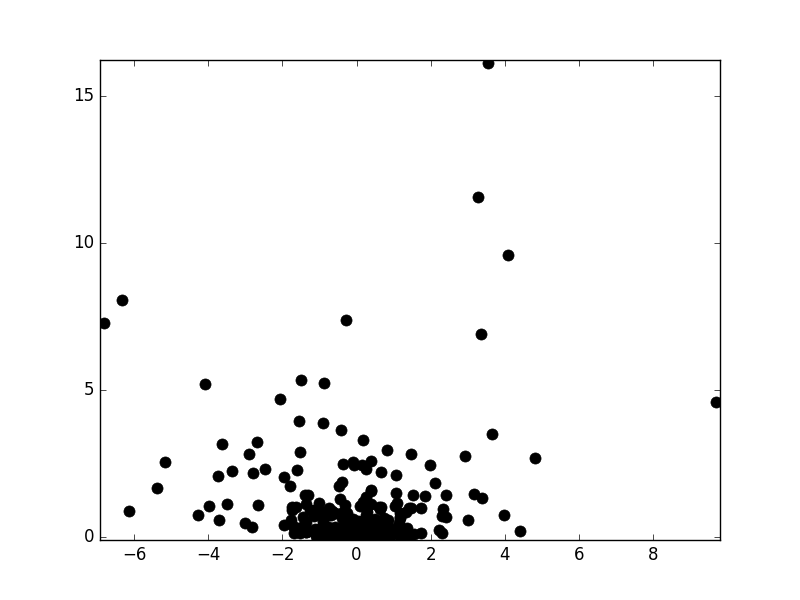}}
		{\includegraphics[width=60mm,height=50mm]{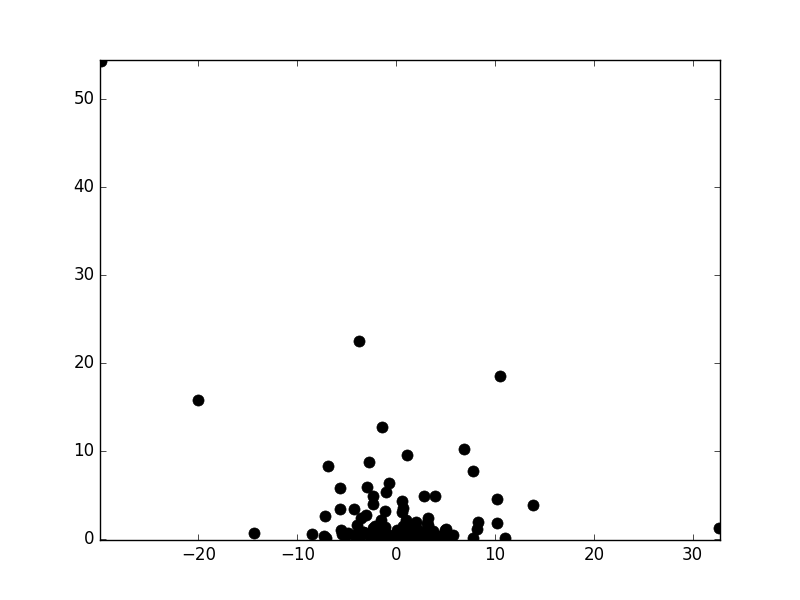}}
		{\includegraphics[width=60mm,height=50mm]{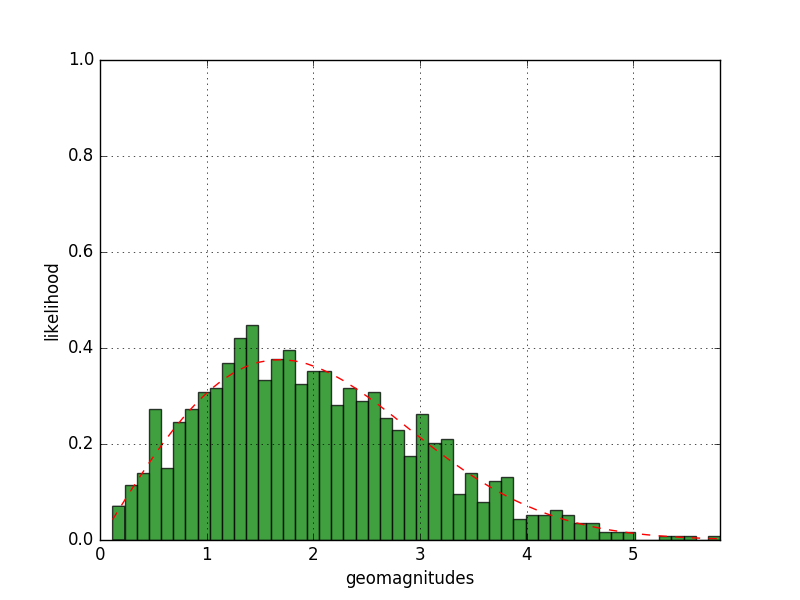}}
		{\includegraphics[width=60mm,height=50mm]{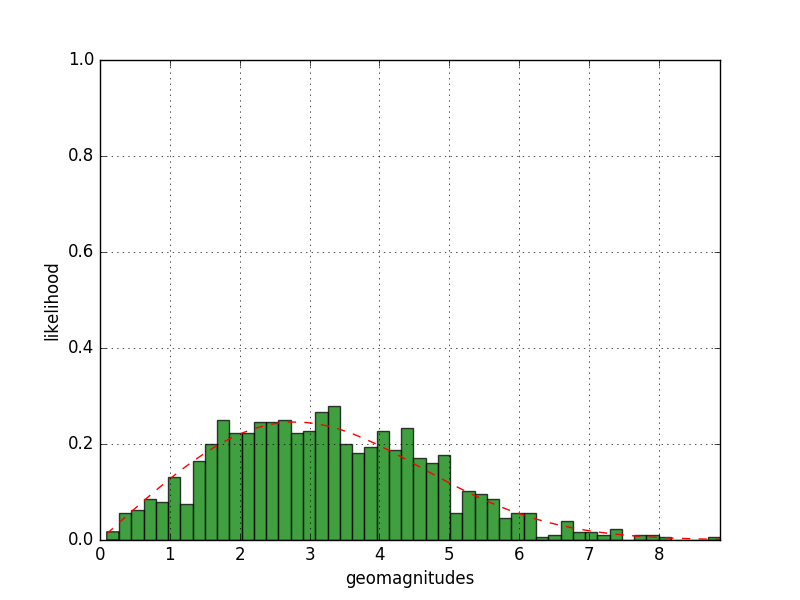}}
  \end{center}
  \caption{Sampling $300$ points from a Gaussian $f_\zeta$ on $\HH$, characterized by $\mathcal{H}[f_\zeta]=2\pi e^{s(s-2\zeta)}$, with parameter $\zeta=1$ (left) and $\zeta=2$ (right). Sampling is implemented via Monte Carlo integration.  The top row illustrates the covariates in $\HH$.  The bottom row illustrates the associated geopolar magnitudes of the covariates.}
  \label{fig:samples}
\end{figure}

A density example is a \technical{Gaussian} $K_{\mathrm{Gauss}}$ defined by
\begin{equation*}
	\hat{K}_{Gauss}(s,b)\propto e^{s(s-2\zeta)},
\end{equation*}
for each choice of real parameter $\zeta$, so named as it is a $\K$-invariant solution to the heat equation (Fig \ref{fig:samples}).  

A $\G$-kernel on $\X$ is \technical{$(\beta,\gamma)$-smooth} if there exist $C_2,C_3>0$ such that
$$C_2e^{-\frac{|s|^{\beta}}{\gamma}}\leqslant |\mathcal{H}K(s,b)|\leqslant C_3 e^{-\frac{|s|^{\beta}}{\gamma}}.$$
for all $s=\rho_{\X}+i\lambda$ and $b\in\B$.
Examples include the previous two $\G$-kernels.

For each $\G$-kernel $K$ on $\X$ and a \technical{bandwidth} parmeter $h>0$, define $K_h$ by
$$(\mathcal{H}K_h)(\lambda,b)=(\mathcal{H}K)(h\lambda,b).$$
The \technical{$\G$-kernel density estimator} $f^{(h)}_{x_1,\ldots,x_n}:\X\rightarrow\R$ is 
$$\hat{f}^{(h)}_{x_1,\ldots,x_n}=\mathcal{H}^{-1}\left[\empchar_{x_1,\ldots,x_n}\mathcal{H}[K_h]\right],$$
for each choice of $\G$-kernel $K$, bandwidth $h>0$, and $x_1,\ldots,x_n\in\X$.

\subsection{Main theorem}
A proof of the following main result can be found in Section \ref{sec:proofs}.

\begin{thm}
  \label{thm:upper}
  For density $f_X\in\mathcal{F}_\alpha(Q)$ and $(\beta,\gamma)$-smooth kernel $K$ on $\X$, 
  \begin{equation*}
		\mathbb{E}\|f_{X_1,\ldots,X_n}^{(h)}-f_X\|^2 \leqslant C_1Q\rho_{\X}^{2\alpha}h^{2\alpha}+K_1QT^{-2\alpha}+K_2n^{-1}e^{-2(h\rho_{\X})^\beta/\gamma}T^{\dim\X},\quad
		X_1,\ldots,X_n\sim_{iid} f_X
  \end{equation*}
	where $C_1,K_1,K_2>0$ are constants not dependent on $\alpha,Q,n$, for each choice of $T>0$.
\end{thm}

By choosing a smooth enough kernel $K$ and optimal bandwidth $h$, we obtain the following upper bound.

\begin{cor}
  For density $f_X\in\mathcal{F}_\alpha(Q)$ and $(\beta,\gamma)$-smooth kernel $K$ on $\X$,
  \begin{equation*}
		\mathbb{E}\|f_{X_1,\ldots,X_n}^{(h)}-f_X\|^2 \leqslant Cn^{-2\alpha/(2\alpha+\dim{\X})},
    \quad
    X_1,\ldots,X_n\sim_{iid} f_X,
  \end{equation*}
  for some constant $C>0$ not dependent on $\alpha,Q,n$ and $h\in\mathcal{O}(n^{-1/(2\alpha+\dim{\X})})$.
  \label{cor:opt}
\end{cor}





\section{Implementation}\label{sec:H2}
For $\X=\HH$ and $\K$ the Gaussian with $\zeta=1$, 
\begin{equation*}
	f^{(h)}_{Z_1,\ldots,Z_n}(z)=\sum_{i=1}^n\int_{-\infty}^{+\infty}\int_{0}^{2\pi}\!\!\mathrm{Im}(k_\theta({Z_i}))^{\frac{1}{2}-i\lambda}e^{-(\frac{h^2}{4}+h^2\lambda^2)}(\mathrm{Im}(k_\theta({z})))^{\frac{1}{2}+i\lambda}\,d\mu,
\end{equation*}
where $d\mu=\frac{1}{8\pi^2n}(\lambda\tanh (\pi\lambda))\,d\theta\,d\lambda$ and $k_\theta$ denotes the rotation matrix associated to the angle $\theta\in[0,2\pi)$.
This double integral, computationally cumbersome, is currently the simplest form known for expressing the KDE - the usual simplifications that allow us to regard a KDE as a convolution of a kernel with an average of dirac point-masses do not work for general $\X$.
We leave for future work the task of optimizing the numerical approximation of the KDE, likely using analogues of discrete Fourier coefficients to construct discrete optimizations of the inverse Helgason-Fourier transform \citep{Pesenson}.

For the present paper, we simulate $1000$ covariates from a Gaussian $f_\zeta$ with dispersion parameter $\zeta=1$ on $\HH$ and numerically compare KDEs defined on the first $n=100,200,\ldots,1000$ covariates with one another and with the true density.
Every $\SO_2$-invariant density is determined by its marginal on geopolar magnitudes.
Given that $f$ is $\SO_2$-invariant, we simplify our computational task by only comparing the associated marginals on geopolar magnitudes.
The associated marginal for $f_\zeta$ is
\begin{equation*}
  2\pi f_\zeta(e^{-r}i)=2\pi(4\pi\zeta)^{-\nicefrac{3}{2}}\sqrt{2}e^{-\zeta/4}\int_{r}^\infty \frac{be^{-b^2/(4\zeta)}}{\sqrt{\cosh(b)-\cosh(r)}}\;db.
\end{equation*}

The associated marginal for the KDE $\hat{f}^h_{X_1,\ldots,X_n}$ is given by
\begin{equation*}
	\hat{f}^h_{X_1,\ldots,X_n}(e^{-r}i)=\sum_{i=1}^n\int_{-\infty}^{+\infty}(\mathrm{Im}Z_i)^{\frac{1}{2}-i\lambda}e^{-(\frac{1}{4}+h^2\lambda^2)}P_{-\half+i\lambda}(\cosh r)\,d\mu.
\end{equation*}
where  $d\mu=\frac{1}{8\pi^2n}(\lambda\tanh (\pi\lambda))\,d\theta\,d\lambda$ and $P_a(c)$ is the Legendre function defined by the integral
\begin{equation*}
  P_a(c)=\frac{1}{2\pi}\int_{0}^{2\pi}(c+\cos\theta\sqrt{c^2-1})^{a}\;d\theta.
\end{equation*}

Bandwidth selection is made by the rule-of-thumb
\begin{equation*}
  h=n^{-\alpha/(\alpha+\dim\X)}
\end{equation*}
instead of a data-driven criterion like cross-validation, due to the computational burden of the formula for the KDE.
An empirical calculation shows that the rule-of-thumb bandwidth closely tracks the ISE minimizer
$$h^*=\arg\min\;\!\!_h\|f-\hat{f}^h_{X_1,\ldots,X_n}\|_2.$$

\begin{figure}
  \begin{center}
    \includegraphics[width=60mm,height=50mm]{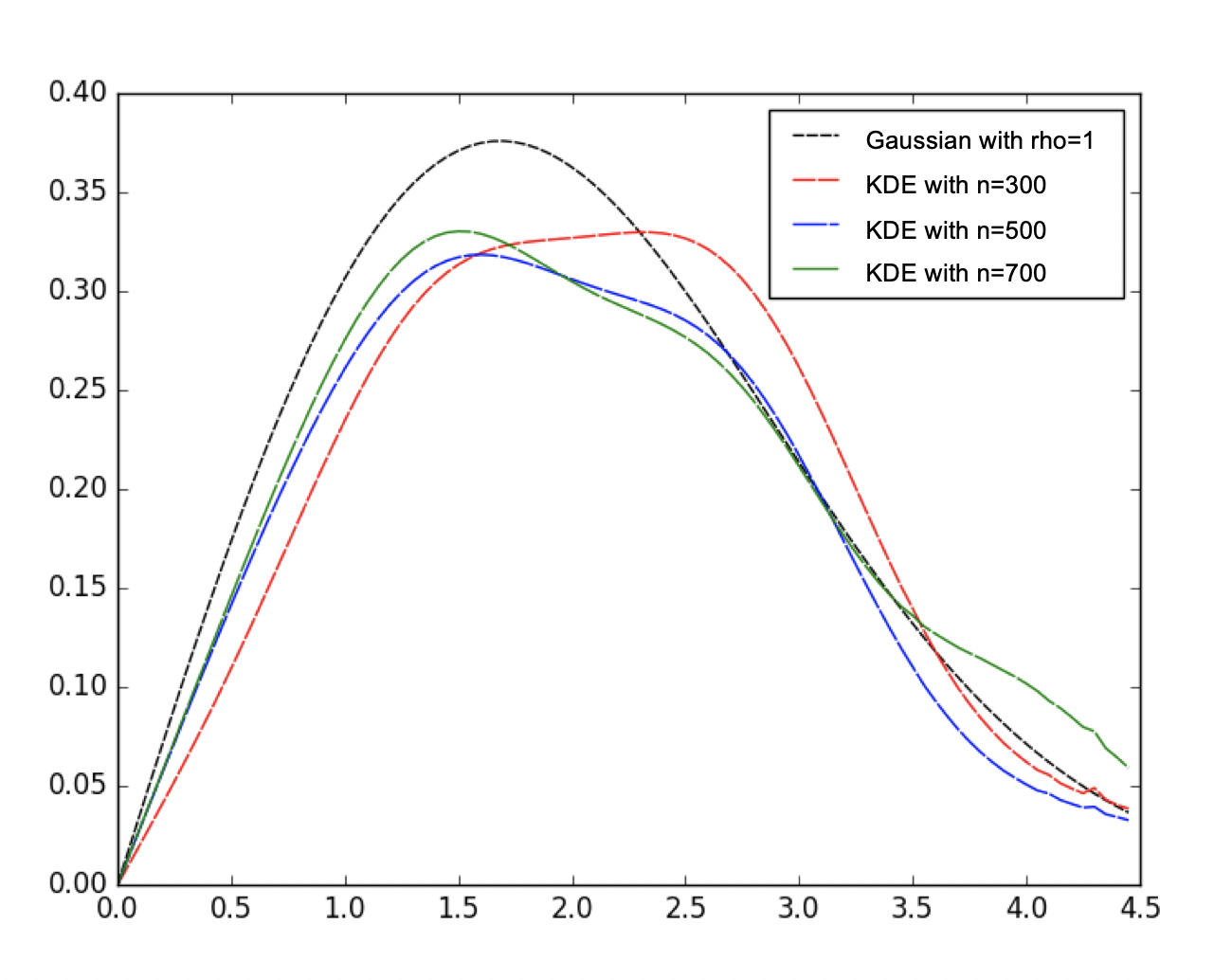}
    \includegraphics[width=60mm,height=50mm]{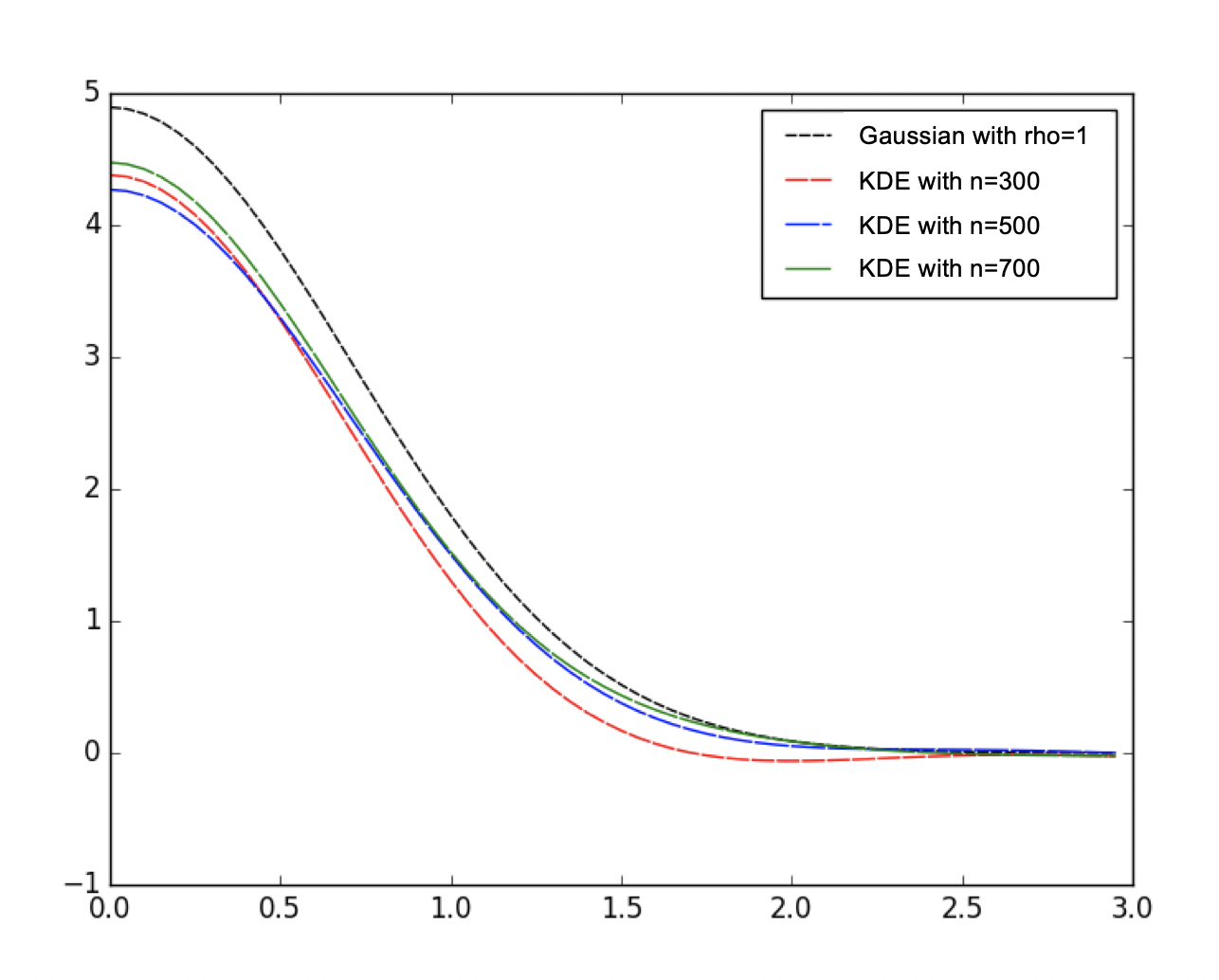}
   \end{center}
	  \caption{Comparing KDEs. The true density of a Gaussian on $\HH$ with dispersion parameter $\zeta=1$ is compared with KDEs, whose bandwidths are chosen by the rule-of-thumb derived in the proof of the minimax rate.  Illustrated above are true and estimated densities with respect to both Lebesgue measure (left) and the measure induced from the volume measure (right), on geopolar magnitudes conditioned on geopolar angles being $0$.  The densities on the right are used to calculate risk, which can be seen to be decreasing in the number of sample points.}
	  \label{fig:KDEs}
\end{figure}

\begin{figure}
  \begin{center}
    \mbox{\includegraphics[width=60mm,height=50mm]{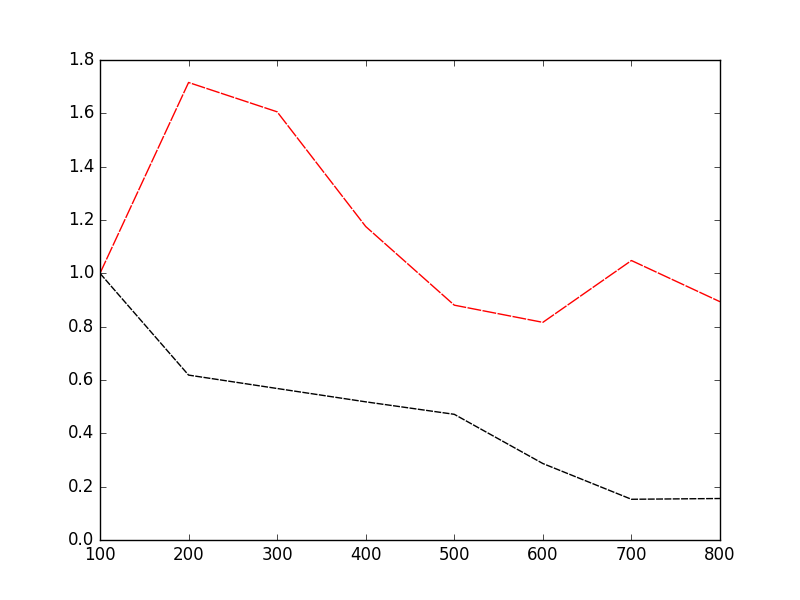}}
   \end{center}
  \caption{Classical versus Non-classical KDEs. The integrated squared error (ISE) with respect to the volume measure on $\HH$ is computed for a KDE on $\HH$ (black) versus a Euclidean KDE (red) for $n=300,500,700$ (x-axis) and normalized so that the ISEs for $n=300$ are both $1$.}
  \label{fig:risk-plots}
\end{figure}

We plot some KDEs at $n=300,500,700$ against the marginal of the true density $f_\zeta$ as functions on $\HH$ (Fig \ref{fig:KDEs}) and as associated Helgason-Fourier transforms on $\R\times\SO_2$ (Fig \ref{fig:KDEs}).
We then compare the rates at which the integrated square errors for the marginals of the KDEs and a Euclidean KDE on sample geopolar magnitudes decrease.  
Our bandwidth-selection criterion for the Euclidean KDE is Silverman's rule-of-thumb.  
We plot the Integrated Squared Error (ISE) for the marginal of the KDE and the Euclidean KDE on the space of geopolar magnitudes (Fig \ref{fig:risk-plots}).

\section{Proofs}\label{sec:proofs}
We use the mean integrated squared error to measure the performance of our generalized estimator.  
As previously noted, proofs here adapt techniques developed in \citep{Huckemann-et-al}.
We break the mean integrated squared error into two parts, variance and squared bias, and bound each part separately.
Throughout, let $s$ denote $\rho_{\X}+i\lambda$, $T$ denote a real number, and $\specmeasure$ denote the spectral measure $|c(\lambda)|^{-2}d\lambda db$.

{\bf Proof of Main Theorem:}
Let $V=\mathbb{E}\left\|f_X^{(n,h)}-\mathbb{E}\left[f_X^{(n,h)}\right]\right\|^2$.  
To obtain the bound on the variance $V$, note that for each $T$,
\begin{align*}
   V
  =\;&\mathbb{E}\int_{\X}\left|f_X^{(n,h)}-\mathbb{E}\left[f_X^{(n,h)}\right]\right|^2dx
  =w_{\X}^{-1}\mathbb{E}\int_{\mathfrak{a}^*}\int_{\B}\left|\mathcal{H}\left[f_X^{(n,h)}\right]-\mathbb{E}\left[\mathcal{H}[f_X^{(n,h)}\right]\right|^2d\specmeasure\\
  =\;&w_{\X}^{-1}\mathbb{E}\int_{\mathfrak{a}^*}\int_{\B}\left|\empchar_{X_1,\ldots,X_n}\mathcal{H}K_h-\mathcal{H}f_X\mathcal{H}K_h\right|^2d\specmeasure
  =w_{\X}^{-1}\mathbb{E}\int_{\mathfrak{a}^*}\int_{\B}\left|\mathcal{H}K_h\right|^2\left|\empchar_{X_1,\ldots,X_n}-\mathcal{H}f_X\right|^2d\specmeasure\\
  =\;&w_{\X}^{-1}\int_{\mathfrak{a}^*}\int_{\B}\left|\mathcal{H}K_h\right|^2\mathbb{E}\left[|\empchar_{X_1,\ldots,X_n}|^2+|\mathcal{H}f_X|^2-2|\empchar||\mathcal{H}f_X|\right]d\specmeasure
\end{align*}

Let $I=\mathbb{E}\left[|\empchar_{X_1,\ldots,X_n}|^2+|\mathcal{H}f_X|^2-2|\empchar||\mathcal{H}f_X|\right]$.
Observe:
\begin{align*}
  I=\;&\left[|\mathcal{H}f_X|^2+n^{-1}\left(\left|\mathbb{E}[e^{2\bar{s}(a(X,b))}]\right|-|\mathcal{H}f_X|^2\right)+|\mathcal{H}f_X|^2-2\mathcal{H}f_X\mathbb{E}[\empchar]\right]\\
	\leqslant\;&\left[|\mathcal{H}f_X|^2+n^{-1}\left(\left|\mathbb{E}[e^{2\rho_{\X}(a(X,b))}]\right|-|\mathcal{H}f_X|^2\right)+|\mathcal{H}f_X|^2-2\mathcal{H}f_X\mathbb{E}[\empchar]\right]\\
  =\;&\left[ 2|\mathcal{H}f_X|^2+n^{-1}\left(|\mathcal{H}f_X(2\rho_{\X},b)|-|\mathcal{H}f_X|^2\right)-2|\mathcal{H}f_X|^2 \right]
  =n^{-1}\left(|\mathcal{H}f_X(2\rho_{\X},b)|-|\mathcal{H}f_X|^2\right)\\
	\leqslant\;&n^{-1}|\mathcal{H}f_X(2\rho_{\X},b)|.
\end{align*}

It then follows that
\begin{align*}
  V
  \leqslant\;&n^{-1}w_{\X}^{-1}\int_{\mathfrak{a}^*}\int_{\B}\left|\mathcal{H}K_h\right|^2|\mathcal{H}f_X(2\rho_{\X},b)|\;d\specmeasure\\
  \leqslant\;&n^{-1}w_{\X}^{-1}\left(\int_{\mathfrak{a}^*}\left|\mathcal{H}K_h\right|^2|c(\lambda)|^{-2}\;d\lambda\right)\left(\int_{\B}|\mathcal{H}f_X(2\rho_{\X},b)|\;db\right).
\end{align*}

The last term is an element inside
\begin{equation*}
	\mathcal{O}\left(\frac{1}{n}\left(\int_{|\lambda|<T}\left|\mathcal{H}K_h\right|^2|c(\lambda)|^{-2}\;d\lambda+\int_{|\lambda|>T}\!\!\!\!\!\!\!\left|\mathcal{H}K_h\right|^2|c(\lambda)|^{-2}\;d\lambda\right)\right),
\end{equation*}
therefore an element inside
\begin{equation*}
	\mathcal{O}\left(\frac{1}{n}\left(e^{-2(h\rho_{\X})^\beta/\gamma}(1+T)^{\dim\mathfrak{n}}\;T^{\dim\mathfrak{a}}+\int_{|\lambda|>T}\!\!\!\!\!\!\!e^{-2h^\beta|\lambda|^\beta/\gamma}(1+|\lambda|)^{\dim\mathfrak{n}}\;d\lambda\right)\right),
\end{equation*}
and hence an element inside
\begin{equation*}
  \mathcal{O}\left(\frac{1}{n}e^{-2(h\rho_{\X})^\beta/\gamma}T^{\dim\X}\right).
\end{equation*}

Let $B=\left\|\mathbb{E}f_{X_1,\ldots,X_n}^{(h)}-f_X\right\|$.  
To obtain the bound on the squared bias $B^2$, note that for each $T$,
\begin{align*}
    B^2
  \;=&\int_{\X}\left|\mathbb{E}f_{X_1,\ldots,X_n}^{(h)}-f_X\right|^2dx
  =w_{\X}^{-1}\int_{\mathfrak{a}^*}\int_{\B}\left|\mathbb{E}\mathcal{H}f_{X_1,\ldots,X_n}^{(h)}-\mathcal{H}f_X\right|^2d\specmeasure
  =w_{\X}^{-1}\int_{\mathfrak{a}^*}\int_{\B}\left|\mathcal{H}f_X\mathcal{H}K_h-\mathcal{H}f_X\right|^2d\specmeasure\\
  \;=&w_{\X}^{-1}\int_{\mathfrak{a}^*}\int_{\B}\left|\mathcal{H}f_X\right|^2\left|\mathcal{H}K_h-1\right|^2d\specmeasure
  =w_{\X}^{-1}\int_{\mathfrak{a}^*}\int_{\B}\overline{s(s-2\rho_{\X})}^{-\alpha}\overline{s(s-2\rho_{\X})}^{\alpha}\left|\mathcal{H}f_X\right|^2\left|\mathcal{H}K_h-1\right|^2d\specmeasure\\
  \;\leqslant&w_{\X}^{-1}\int_{|\lambda'|<hT}\int_{\B}h^{2\alpha-1}(-(\lambda')^2-\rho_{\X}^2)^{-\alpha}\left|\mathcal{H}K-1\right|^2 d\specmeasure'
  +\mathcal{O}\left(\sup_{|\lambda|>T}|s|^{-2\alpha}\left|\mathcal{H}K_h-1\right|^2\right)\\[.1in]
  \;\subset&\mathcal{O}\left(C_1h^{2\alpha}+K_1T^{-2\alpha}\right)
\end{align*}
for constants $C_1,K_1>0$, where $\lambda'=h\lambda$ and $\specmeasure'=h\specmeasure$ above.
\vspace{.1in}\\
{\bf Proof of the Corollary:}
By Theorem \ref{thm:upper}, for each $T>0$
\begin{align*}
  \mathbb{E}\|f_X^{(n,h)}-f_X\|^2 &\leqslant C_1h^{2\alpha}+K_1T^{-2\alpha}+K_2n^{-1}e^{-2(h\rho_{\X})^\beta/\gamma}T^{\dim\X}.
\end{align*}
Setting $T^{-2\alpha}\propto n^{-1}e^{-2(h\rho_{\X})^\beta/\gamma}T^{\dim\X}$, we obtain
$$T\propto n^{1/(2\alpha+\dim\X)}.$$

The upper bound converges at the fastest possible rate when
$$h^{2\alpha}\propto T^{-2\alpha}\propto n^{-2\alpha/(2\alpha+\dim\X)},$$
or equivalently, when $h\propto n^{-\alpha/(\alpha+\dim\X)}$.  
Thus
\begin{align*}
  \mathbb{E}\|f_X^{(n,h)}-f_X\|^2\in\mathcal{O}(n^{-2\alpha/(2\alpha+\dim\X)}).
\end{align*}


\section{Conclusion}
Until now, kernel density estimation in the non-compact setting has required either that one restricts to Euclidean space, requires that the densities have compact support, or requires that one specifies a kernel-like function for each point in a complete Riemannian manifold.
We have introduced a new density estimator on a large class of non-compact symmetric spaces, sidestepping these various restrictions, and have proven an upper bound for the rate of convergence, conjecturally a minimax rate, identical to the minimax rate of convergence for a Euclidean kernel density estimator. 
We then specialized our generalized kernel density estimator for the hyperboloid, motivated by applications to network inference. 
Future work will explore adaptivity, computational optimizations, and applications to symmetric spaces other than the hyperboloid.

\section*{Acknowledgements}
This work was partially supported by an NSF Graduate Research Fellowship under grant DGE-1252522. Also, the author is grateful to her advisor, Cosma Shalizi, for his invaluable guidance and feedback throughout this research. 
The author is also grateful to the anonymous referee for helpful comments, corrections, and suggestions.  

\bibliographystyle{abbrvnat}
\bibliography{bib_asta}

\end{document}